\documentclass{ifacconf}

\usepackage{amsmath}        % required for math environments
\usepackage{graphicx}       % required for figures
\usepackage{natbib}         % required for bibliography
\usepackage{epstopdf}       % required for eps figures

\begin{document}
    \begin{frontmatter}
        % DOCUMENT INFORMATION
        \title{Large-scale Virtual Clinical Trials\\of Closed-loop Treatments\\for People with Type 1 Diabetes}

\thanks[footnoteinfo]{This work was partially funded by the IFD Grand Solution project
ADAPT-T2D (9068-00056B).}

\author{Tobias K. S. Ritschel,}
\author{Asbj{\o}rn Thode Reenberg,}
\author{John Bagterp J{\o}rgensen}

\address{Department of Applied Mathematics and Computer Science,\\Technical University of Denmark, DK-2800 Kgs. Lyngby, Denmark (e-mail: jbjo@dtu.dk).}

        % ABSTRACT
        \begin{abstract}
            We propose a virtual clinical trial for assessing the safety and efficacy of closed-loop diabetes treatments prior to an actual clinical trial. Such virtual trials enable rapid and risk-free pretrial testing of algorithms, and they can be used to compare different treatment variations for large and diverse populations. The participants are represented by multiple mathematical models, consisting of stochastic differential equations, and we use Monte Carlo closed-loop simulations to compute detailed statistics of the closed-loop treatments. We implement the virtual clinical trial using high-performance software and hardware, and we present an example trial with two mathematical models of one~million participants over 52~weeks (i.e., two~million simulations), which can be completed in 2~h 9~min.

        \end{abstract}
        
        % KEYWORDS
        \begin{keyword}
            Virtual clinical trials \sep
Diabetes \sep
Stochastic modeling \sep
High-performance computing

        \end{keyword}
    \end{frontmatter}

    % MAIN CONTENT
    \section{Introduction}\label{sec:introduction}
Clinical trials ensure the safety and efficacy of medical treatments, but they are also time-consuming and expensive. Furthermore, they might result in a negative outcome, e.g., that the proposed treatment is not safe.
Therefore, prior to the actual clinical trial, it is important to 1)~evaluate the potential treatment performance, 2)~identify shortcomings and risks, and 3)~assess the advantages of the treatment over alternative treatments.
This is the purpose of \emph{virtual} clinical trials which involve virtual participants that are represented by a mathematical model (often consisting of differential equations). By using high-performance computing (HPC) software and hardware, such virtual trials can involve large populations and long-term protocols, which allows for extensive and fast testing of different treatment variations.

In this paper, we consider virtual clinical trials of closed-loop diabetes treatment systems. These are also referred to as \emph{artificial pancreases} (APs). Worldwide, one in ten adults suffer from diabetes, and according to the~\citet{IDF:2021}, it accounted for 9\% of the 2021 global health expenditure (USD 966~billion). Specifically, we consider type 1 diabetes (T1D) where, due to autoimmune destruction of the $\beta$-cells, the pancreas does not produce any insulin. People with T1D require daily treatment with exogenous insulin in order to avoid high blood glucose concentrations (hyperglycemia). Prolonged hyperglycemia can lead to a number of health complications and chronic conditions, e.g., chronic kidney disease, cardiovascular disease, and damage to the eyes and nerves. Conversely, the insulin treatment can lead to low blood glucose concentrations (hypoglycemia), which can result in acute complications such as loss of consciousness and seizures.

Clearly, given the risks associated with hyper- and hypoglycemia, it is not straightforward for people with T1D to manage their insulin treatment. Therefore, over the last few decades, there have been significant developments within AP systems which can decrease this burden~\citep{Lal:etal:2019}. APs typically consist of 1)~a continuous glucose monitor (CGM), 2)~an insulin pump, and 3)~a control algorithm implemented on a smartphone or a dedicated device. The control algorithm repeatedly computes an appropriate insulin flow rate based on measurements from the CGM device and communicates it to the insulin pump. There exist a variety of control algorithms for computing the insulin flow rate, e.g., based on heuristics, fuzzy logic, proportional-integral-derivative (PID) control~\citep{Huyett:etal:2015, Sejersen:etal:2021}, and model predictive control (MPC)~\citep{Boiroux2:etal:2018, Chakrabarty:etal:2020}.
All of these algorithms contain algorithmic parameters which must be tuned based on simulation, i.e., based on a virtual clinical trial. As the human physiology and behavior vary significantly between people and over time, this is a nontrivial task. In spite of this, the tuning is typically based on short-term simulations of one or a few virtual participants who are only represented by a single mathematical model. In contrast, if large-scale long-term virtual clinical trials (involving multiple mathematical models) are used to identify candidate algorithms and algorithmic parameters, the chances of a successful real-world clinical trial increase significantly.

In this work, we develop an approach for performing large-scale long-term virtual clinical trials of AP algorithms (i.e., Monte Carlo closed-loop simulations). The approach is an extension of our recent work~\citep{Reenberg:etal:2022, Wahlgreen:etal:2021} which allows the virtual participants to be represented by multiple mathematical models. Specifically, we use 1)~the model described by~\citet{Hovorka:etal:2002} combined with a CGM model and 2)~a modification of the model described by~\citet{DallaMan:etal:2007} and~\citet{Colmegna:etal:2020}. Furthermore, we extend both models to incorporate uncertainty about the physiology and the CGM measurements, i.e., we formulate the models as stochastic differential equations (describing the dynamics) and stochastic algebraic equations (describing the CGM measurements). We implement the approach using high-performance C code, and we use an HPC cluster~\citep{DTU_DCC_resource} to carry out the computations. Finally, we demonstrate the utility of the virtual clinical trial by comparing the performance of an AP algorithm for the two models. Specifically, we simulate one~million participants over 52~weeks with both models (i.e., two~million in total) which takes 2~h 9~min.

The remainder of this paper is structured as follows. In Section~\ref{sec:virtual:clinical:trial}, we describe the virtual clinical trial, and in Section~\ref{sec:models}, we describe the two models used in this work. Next, we present the results of the virtual clinical trial in Section~\ref{sec:results}, and finally, we present conclusions in Section~\ref{sec:conclusions}.

    \section{Virtual clinical trial}\label{sec:virtual:clinical:trial}
The virtual clinical trial involves a population of people with T1D, a protocol describing, e.g., the size and duration of meals, mathematical models, model parameter values, and AP algorithms. Additionally, meals and other activities can be incorrectly announced to the AP algorithm (as is often the case in reality). Furthermore, the mathematical models can be either deterministic (no uncertainty) or stochastic (uncertain dynamics and measurements). The uncertainty can, for instance, represent physiological phenomena that are not included in the model or unknown model parameters.

\subsection{Virtual people}\label{sec:virtual:clinical:trial:people}
The virtual clinical trial contains one~million virtual people. Each person is represented by the same pieces of information as real people, e.g., a unique ID, given name, family name, date of birth, place of birth, sex, height, and body weight. Additionally, each person is associated with a set of mathematical models and model parameter values.

\subsection{Protocols}\label{sec:virtual:clinical:trial:protocols}
A protocol describes the participants' activities during the virtual trial. Typical examples include the time and carbohydrate contents of meals as well as the duration and intensity of physical activity. Furthermore, each protocol has its own ID, and it contains IDs, time stamps (start and end time), and type and size for each activity.

In Section~\ref{sec:results}, we demonstrate the virtual clinical trial using a protocol based on a Northern European lifestyle with respect to 1)~meal times, 2)~work weeks, 3)~vacation weeks, 4)~public holidays, and 5)~seasons. The year is divided into 4 seasons consisting of 13 weeks. In total, there are 8~weeks of vacation (2~of them represent public holidays). All weeks are either a standard week, an active week, or a vacation week. Similarly, all days are categorized as either a standard day, an active day, a day with a movie night, or a day with a late night. Each type of week contains a different combination of the days, and each season contains a different combination of the weeks. Table~\ref{tab:WeeksAndSeasons} gives an overview of the weeks and seasons as well as the carbohydrate contents of the different meals in the trial. During autumn, winter, and vacation weeks, the participants are less active and eat more. Furthermore, in addition to the meals consumed during a standard day, active days have an additional exercise session, days with a movie night have an additional snack in the evening, and days with a late night have two additional snacks in the evening.
\begin{table}
	\setlength{\tabcolsep}{3pt}
	\caption{The seasons, weeks, and meal carbohydrate contents~\citep{Reenberg:etal:2022}.}
	\label{tab:WeeksAndSeasons}
	\begin{tabular}{p{0.15\linewidth}|ccc}
		\multicolumn{4}{l}{\textbf{Compositions of the seasons}} \\[2pt]
		\hline
		Season & Standard weeks & Active weeks & Vacation weeks \\
		\hline
		Winter & 6 & 4 & 3 \\
		Spring & 6 & 6 & 1 \\
		Summer & 7 & 3 & 3 \\
		Autumn & 9 & 3 & 1 \\
		\hline
	\end{tabular}

	\vspace{7pt}

	\begin{tabular}{p{0.16\linewidth}|cccc}
		\multicolumn{5}{l}{\textbf{Compositions of the weeks}} \\[2pt]
		\hline
		Week type & Standard day & Active day & Movie nights & Late nights \\
		\hline
		Standard & 4 & 1 & 1 & 1  \\
		Active   & 3 & 3 & 1 & 0  \\
		Vacation & 5 & 0 & 0 & 2  \\
		\hline
	\end{tabular}

	\vspace{7pt}

	\begin{tabular}{l|cc}
		\multicolumn{3}{l}{\textbf{Body weight-dependent meal carbohydrate contents}} \\[2pt]
		\hline
		Meal size 	& Amount of carbohydrates & For a 70~kg person \\
		\hline
		Large meal  & 1.29~g CHO/kg & 90 g CHO \\
		Medium meal & 0.86~g CHO/kg & 60 g CHO \\
		Small meal  & 0.57~g CHO/kg & 40 g CHO \\
		Snack       & 0.29~g CHO/kg & 20 g CHO \\
		\hline
	\end{tabular}
\end{table}
\begin{figure}[tb]
	\centering
	\includegraphics[width=0.45\textwidth]{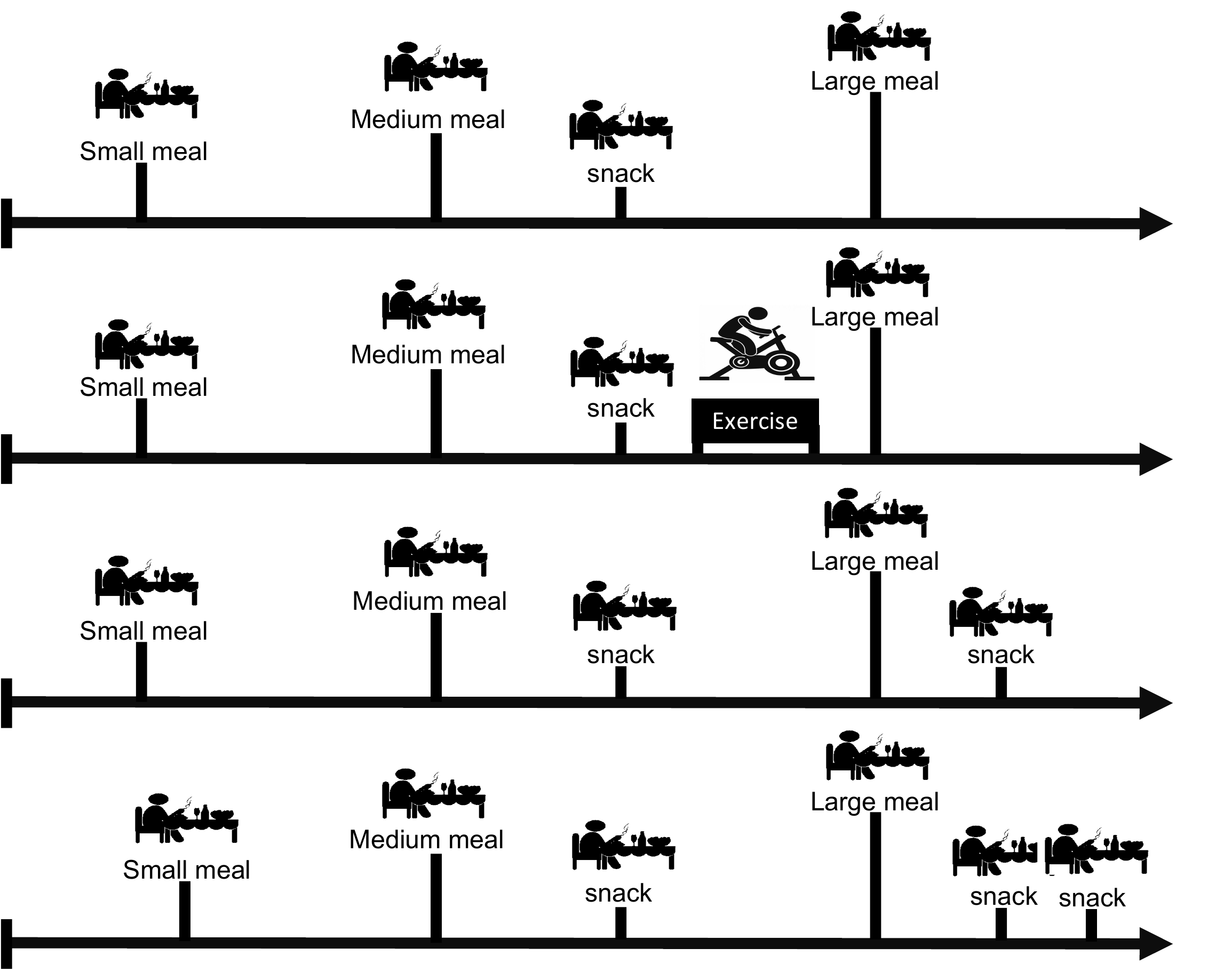}
	\caption{Schematic of the different types of days during autumn and winter~\citep{Reenberg:etal:2022}. Top: A standard day. Second from the top: An active day. Third from the top: A day with a movie night. Bottom: A day with a late night. In the spring and summer, the snack is consumed between breakfast and lunch, and the dinner is a medium meal.}
	\label{fig:WinterSummerDays}
\end{figure}

\subsection{Database}
The virtual participants, model parameters, protocols, and results from the virtual clinical trial are stored in a database. We use PostgreSQL which is an open-source database system. The database enables direct comparison of the performance of different AP systems on different populations, e.g., people with high body weights or low insulin sensitivity (provided that this sensitivity is a model parameter). Additionally, the database contains basic components for building protocols, e.g., the 4~types of days used in the protocol described in Section~\ref{sec:virtual:clinical:trial:protocols}. The user can construct their own protocols based on these basic components. Finally, the database can be used together with a graphical user interface to select and visualize the results of the virtual clinical trial, show statistics of certain demographics and protocols, and add new elements (for instance, virtual people, mathematical models and parameters, or protocols).

    \section{Models}\label{sec:models}
In this section, we describe the two models used in the virtual clinical trial. In Section~\ref{sec:models:hovorka}, we describe the model developed by~\citet{Hovorka:etal:2002} combined with a CGM model, and in Section~\ref{sec:models:uva:padova}, we describe the model developed by~\citet{DallaMan:etal:2007} and~\citet{Colmegna:etal:2020} where the meal model is replaced with that proposed by~\citet{Hovorka:etal:2002}.
In order to simplify the demonstration of the virtual clinical trial in Section~\ref{sec:results}, we do not include exercise in the models or in the protocol.

\subsection{An extension of Hovorka's model}\label{sec:models:hovorka}
The insulin subsystem is described by
\begin{subequations}\label{eq:hovorka:insulin}
	\begin{align}
		\label{eq:hovorka:insulin:s1}
		\dot S_1(t) &= u_I(t)    - \frac{S_1(t)}{\tau_S}, \\
		\label{eq:hovorka:insulin:s2}
		\dot S_2(t) &= \frac{S_1(t)}{\tau_S} - \frac{S_2(t)}{\tau_S}, \\
		\label{eq:hovorka:insulin:i}
		\dot I(t)   &= \frac{1}{V_I}\frac{S_2(t)}{\tau_S} - k_e I(t),
	\end{align}
\end{subequations}
where $S_1$ and $S_2$~[mU] constitute a two-compartment chain representing the absorption of insulin, $I$~[mU/L] is the insulin concentration in the plasma, $\tau_S$~[min] is the insulin absorption time constant, $V_I$ is the volume in which insulin is distributed, $k_e$~[1/min] is an elimination rate, and $u_I$~[mU/min] is the insulin infusion rate.
The insulin action subsystem is described the three compartments
\begin{subequations}\label{eq:hovorka:insulin:action}
	\begin{align}
		\label{eq:hovorka:insulin:action:x1}
		\dot x_1(t) &= k_{b1} I(t) - k_{a1} x_1(t), \\
		\label{eq:hovorka:insulin:action:x2}
		\dot x_2(t) &= k_{b2} I(t) - k_{a2} x_2(t), \\
		\label{eq:hovorka:insulin:action:x3}
		\dot x_3(t) &= k_{b3} I(t) - k_{a3} x_3(t),
	\end{align}
\end{subequations}
where $x_i$~[1/min] are the insulin effect on the glucose distribution ($i=1$), the disposal of glucose ($i=2$), and the endogenous glucose production ($i=3$). Furthermore, $k_{bi}$~[(L/mU)/min$^2$] are activation rates and $k_{ai}$~[1/min] are deactivation rates (for $i=1, 2, 3$).
The meal subsystem is described by a two-compartment chain:
\begin{subequations}\label{eq:hovorka:meal}
	\begin{align}
		\label{eq:hovorka:meal:d1}
		\dot D_1(t) &= A_G D(t)    - \frac{D_1(t)}{\tau_D}, \\
		\label{eq:hovorka:meal:d2}
		\dot D_2(t) &= \frac{D_1(t)}{\tau_D} - \frac{D_2(t)}{\tau_D}.
	\end{align}
\end{subequations}
Here, $D$~[mmol/min] is the meal carbohydrate content (per minute), $D_1$ and $D_2$~[mmol] represent the meal absorption, $A_G$~[--] is the bioavailibility of the carbohydrates, and $\tau_D$~[min] is a time constant.
The glucose subsystem is also described by two compartments, i.e.,
\begin{subequations}\label{eq:hovorka:glucose}
	\begin{align}
		\label{eq:hovorka:glucose:q1}
		\dot Q_1(t)
		&= \frac{D_2(t)}{\tau_D} - F_{01,c}(t) - F_R(t) - x_1(t) Q_1(t) \nonumber \\
		&\quad + k_{12} Q_2(t) + EGP_0 (1 - x_3(t))), \\
		\label{eq:hovorka:glucose:q2}
		\dot Q_2(t) &= x_1(t) Q_1(t) - k_{12} Q_2(t) - x_2 Q_2(t),
	\end{align}
\end{subequations}
where $Q_1$ and $Q_2$~[mmol] are the accessible and non-accessible glucose compartments, $k_{12}$~[1/min] is a transfer rate, $EGP_0$~[mmol] is the endogenous glucose production (extrapolated to an insulin concentration of zero). Furthermore,
\begin{subequations}\label{eq:hovorka:glucose:piecewise}
	\begin{align}
		\label{eq:hovorka:glucose:piecewise:f01c}
		F_{01,c}(t) &=
		\begin{cases}
			F_{01} & G(t) \geq 4.5~\mathrm{mmol/L}, \\
			F_{01} G(t)/4.5 & \text{otherwise},
		\end{cases} \\
		\label{eq:hovorka:glucose:piecewise:fr}
		F_R(t) &=
		\begin{cases}
			0.003(G(t) - 9) V_G & G(t) \geq 9~\mathrm{mmol/L}, \\
			0 & \text{otherwise},
		\end{cases} \\
		G(t) &= \frac{Q_1(t)}{V_G},
	\end{align}
\end{subequations}
where $F_{01,c}$ and $F_{01}$~[mmol/min] are the corrected and nominal total non-insulin dependent glucose fluxes, $F_R$ [mmol/min] is the renal glucose clearance, $G$~[mmol/L] is the glucose plasma concentration, and $V_G$~[L] is the volume in which the glucose is distributed.
Finally, the CGM subsystem is
\begin{equation}\label{eq:hovorka:cgm}
	\dot G_I(t) = \frac{G(t)}{\tau_{IG}} - \frac{G_I(t)}{\tau_{IG}}.
\end{equation}
Here, $G_I$~[mmol/L] is the interstitial glucose concentration measured by the CGM and $\tau_{IG}$~[1/min] is a time constant.

\subsection{A modification of the UVA/Padova model}\label{sec:models:uva:padova}
In the UVA/Padova model, the glucose subsystem is described by
\begin{subequations}\label{eq:uva:padova:glucose}
	\begin{align}
	    \label{eq:uva:padova:glucose:gp}
		\dot G_{p} 	&= EGP(t)+Ra_m(t)-U_{ii}(t)-E(t) \nonumber \\ & \quad -k_1G_p(t)+k_2G_t(t), \\
		\label{eq:uva:padova:glucose:gt}
		\dot G_t(t) &=-U_{id}(t)+k_1G_p(t)-k_2G_t(t), \\
		\label{eq:uva:padova:glucose:g}
		G(t) 	    &= \frac{G_p(t)}{V_g},
	\end{align}
\end{subequations}
where $G_p$ and $G_t$~[mg/kg] are the plasma glucose in rapidly and slowly equilibrating tissues, respectively, $EGP$~[(mg/kg)/min] is the endogenous glucose production, $Ra_m$~[(mg/kg)/min] is the glucose rate of appearance, $U_{ii}$ and $U_{id}$~[(mg/kg)/min] are the insulin-independent and insulin-dependent glucose utilization, $E$~[(mg/kg)/min] represents renal excretion, $k_1$ and $k_2$~[1/min] are rate parameters, $G$~[mg/dL] is the glucose concentration, and $V_g$~[dL/kg] is the glucose distribution volume.
Next, the insulin subsystem is represented by
\begin{subequations}
	\begin{align}
		\dot I_{\ell}(t)&= -(m_1+m_3)I_{\ell}(t)+m_2I_p(t), \\
		\dot I_{p}(t) 	&= -(m_2+m_4)I_p(t)+m_1I_\ell(t)+Ra_{Isc}(t),\\
		I(t) &= \frac{I_p(t)}{V_i},
	\end{align}
\end{subequations}
where $I_\ell$ and $I_p$~[pmol/kg] are insulin in the plasma and liver, $m_i$~[1/min] for $i = 1, \ldots, 4$ are rate parameters, $Ra_{Isc}$~[(pmol/kg)/min] is the insulin rate of appearance, $I$~[pmol/L] is the plasma insulin concentration, and $V_i$~[L/kg] is the insulin distribution volume. Furthermore,
\begin{subequations}
	\begin{align}
		m_2 &= \frac{3CL}{5 HE_b V_i BW}, \\
		m_3 &= \frac{HE_b m_1}{1 - HE_b},\\
		m_4 &= \frac{2 CL}{5 V_i BW}.
	\end{align}
\end{subequations}
Here, $CL$~[L/min] is the insulin clearance, $HE_b$~[--] is the basal hepatic insulin extraction, and $BW$~[kg] is the body weight.
The insulin-independent and insulin-dependent glucose utilization are given by
\begin{subequations}
	\begin{align}
		U_{ii}(t)	&= F_{cns}, \\
		U_{id}(t)	&= \frac{(V_{m0} + V_{mx} X(t)) G_t(t)}{K_{m0} + G_t(t)},\\
		\dot X(t) 	&= -p_{2U} X(t) + p_{2U} (I(t) - I_b),
	\end{align}
\end{subequations}
where $F_{cns}$~[(mg/kg)/min] is the glucose uptake of the erythrocytes and the brain, $V_{mx}$~[mg~L/(kg~pmol~min)] and $K_{m0}$~[mg/kg] are parameters, $X$~[pmol/L] is the insulin concentration in the interstitial fluid, $p_{2U}$~[1/min] is the rate of the insulin action on the peripheral glucose utilization, $I_b$~[pmol/L] is the basal insulin plasma concentration, and $V_{m0}$~[(mg/kg)/min] is
\begin{subequations}
	\begin{align}
		V_{m0} &= \frac{(EGP_b - F_{cns})(K_{m0} + G_{tb})}{G_{tb}}, \\
		G_{tb} &= \frac{F_{cns} - EGP_b + k_1 G_{pb}}{k_2},
	\end{align}
\end{subequations}
where $EGP_b$~[(mg/kg)/min], $G_{tb}$~[mg/kg], and $G_{pb}$~[mg/kg] are the basal endogenous glucose production and the basal plasma glucose masses.
The model of the oral glucose absorption is similar to that in~\eqref{eq:hovorka:meal}:
\begin{subequations}
	\begin{align}
		\dot D_1(t) &= A_G D(t) - \frac{D_1(t)}{\tau_D}, \\
		\dot D_2(t) &= \frac{D_1(t)}{\tau_D} - \frac{D_2(t)}{\tau_D},\\
		Ra_m(t) &= \frac{D_2(t)}{BW\tau_D}.
	\end{align}
\end{subequations}
Again, $D$~[mg/min] is the meal carbohydrate rate, $D_1$ and $D_2$~[mg] describe the meal absorption, $A_G$~[--] is the carbohydrate bioavailibility, and $\tau_D$~[min] is a time constant.
The endogenous glucose production is
\begin{subequations}
	\begin{align}
		EGP(t)  		&= \max\{0, EGP_b - k_{p2}(G_p(t) - G_{pb}) \nonumber \\
		&\quad - k_{p3}(I_d(t) - I_b)\},\\
		\dot I_d(t) 	&= -k_i(I_d(t) - I_1(t)),\\
		\dot I_1(t) 	&= -k_i(I_1(t) - I(t)),
	\end{align}
\end{subequations}
where $k_{p2}$~[1/min] and $k_{p3}$~[mg~L/(kg~pmol~min] are the liver glucose effectiveness and the amplitude of the insulin action of the liver, $I_d$ and $I_1$~[pmol/L] constitute a two-compartment delayed insulin signal chain, and $k_i$~[1/min] is a rate parameter.
Next, the renal excretion is given by
\begin{equation}
	E(t) = \max\{0, k_{e1} (G_p(t) - k_{e2})\}.
\end{equation}
Here, $k_{e1}$~[1/min] and $k_{e2}$~[mg/kg] are the glomerular filtration rate and the renal glucose threshold.
The subcutaneous insulin delivery is described by
\begin{subequations}
	\begin{align}
		\dot I_{sc1}(t) &= -k_d I_{sc1}(t) + k_{a1} I_{sc1}(t) + \frac{u_I(t)}{BW}, \\
		\dot I_{sc2}(t) &= k_d I_{sc1}(t) - k_{a2} I_{sc2}(t), \\
		Ra_{Isc}(t) 	&= k_{a1} I_{sc1}(t) + k_{a2} I_{sc2}(t),
	\end{align}
\end{subequations}
where $I_{sc1}$ and $I_{sc2}$~[pmol/kg] are insulin in a non-monomeric and monomeric state, $k_d$, $k_{a1}$, and $k_{a2}$~[1/min] are rate parameters accounting for subcutaneous insulin kinetics, and $u_I$~[pmol/min] is the insulin infusion rate.
Finally, the subcutaneous glucose concentration, $G_{sc}$ [mg/dL], is
\begin{align}
		\dot G_{sc}(t) = -k_{sc} G_{sc}(t) + k_{sc} G(t),
\end{align}
and $k_{sc}$~[1/min] is the inverse of a time constant.

\subsection{General mathematical form and simulation}
The two models described in this section are in the form
\begin{subequations}\label{eq:sde}
	\begin{align}
		\label{eq:sde:initial}
		x(t_0) &= x_0, \\
		\label{eq:sde:state}
		d x(t) &= f(t, x(t), u(t), d(t), p) dt \nonumber \\
		& \quad + \sigma(t, x(t), u(t), d(t), p) dw(t), \\
		\label{eq:sde:output}
		z(t) &= h(t, x(t), p), \\
		\label{eq:sde:observed}
		y(t_k) &= g(t_k, x(t_k), p) + v(t_k).
	\end{align}
\end{subequations}
Here, $t$ is time, $t_0$ is the initial time, $x$ are the states (i.e., the physiological state), $x_0$ are the initial states, $u$ are manipulated inputs computed by the AP algorithm (e.g., the insulin infusion rate), $d$ are disturbance variables (e.g., the meal carbohydrate content), and $p$ are model parameters (specific to each person). The first term in~\eqref{eq:sde:state} is the deterministic (drift) term, and the second term is the stochastic (diffusion) term.
At time $t_k$ (e.g., every 5~min), the AP receives a CGM measurement, $y$, which is corrupted by noise, $v$. Furthermore, the outputs, $z$, are used to evaluate the AP algorithm. The measurement noise is normally distributed, and $w$ is a standard Wiener process, i.e., $v(t_k)\sim N(0, R)$ and $dw(t)\sim N(0, I dt)$ where $I$ is an identity matrix.
In between measurements, the manipulated inputs and the disturbance variables are modeled as constant:
\begin{subequations}
	\begin{align}
		u(t) &= u_k, & t &\in[t_k, t_{k+1}[, \\
		d(t) &= d_k, & t &\in[t_k, t_{k+1}[.
	\end{align}
\end{subequations}
Finally, when the AP algorithm (represented by the functions $\kappa_k$ and $\lambda_k$) receives a measurement, it updates its internal state, $x_k^c$, and computes the manipulated inputs based on the measurement, $y_k = y(t_k)$, the targets $\bar u_k$ and $\bar y_k$, an estimate of the disturbance variables, $\hat d_k$, and the hyperparameters $p_\kappa$ and $p_\mu$:
\begin{subequations}
	\begin{align}
		x_{k+1}^c 	&= \kappa_k(x_k^c, y_k, \bar u_k, \bar y_k, \hat d_k, p_\kappa), \\
		u_k			&= \lambda_k(x_k^c, y_k, \bar u_k, \bar y_k, \hat d_k, p_\mu).
	\end{align}
\end{subequations}

In Section~\ref{sec:models:hovorka} and~\ref{sec:models:uva:padova}, we described the models without uncertainty, i.e., with $\sigma = 0$. In Section~\ref{sec:results}, we use uncertain models where we 1)~add a stochastic term to the plasma glucose compartments~\eqref{eq:hovorka:glucose:q1} and~\eqref{eq:uva:padova:glucose:gp} with (equivalent) diffusion coefficients of 1.5~mmol/min$^{3/2}$ and 270.24/$BW$~mg/(kg~min$^{3/2}$), and 2)~add measurement noise with a variance of 0.1~mmol/L $\approx$ 1.8~mg/dL.

We use the Euler-Maruyama method with a step size of 0.5~min to simulate each participant, and we only store statistics (and not individual simulations). We implement the virtual clinical trial using parallel high-performance C code for shared-memory architectures and two 2.9 GHz AMD EPYC 7542 32-core processors.

    \section{Virtual clinical trial results}\label{sec:results}
In this section, we demonstrate how the virtual clinical trial can be used to assess the safety and efficacy of an AP algorithm. We use the protocol described in Section~\ref{sec:virtual:clinical:trial}, and for each model, we include one million participants. In the analysis, we consider the first 4~weeks a titration period and disregard them. In the following, we refer to the model presented in Section~\ref{sec:models:hovorka} as model~A and the one presented in Section~\ref{sec:models:uva:padova} as model~B. We choose the algorithmic hyperparameters based on deterministic simulations with model~A. The total computation time is 2~h 9~min.

\begin{figure}[tb]
	\centering
	\includegraphics[trim=10 25 85 45, clip, width=0.45\textwidth]{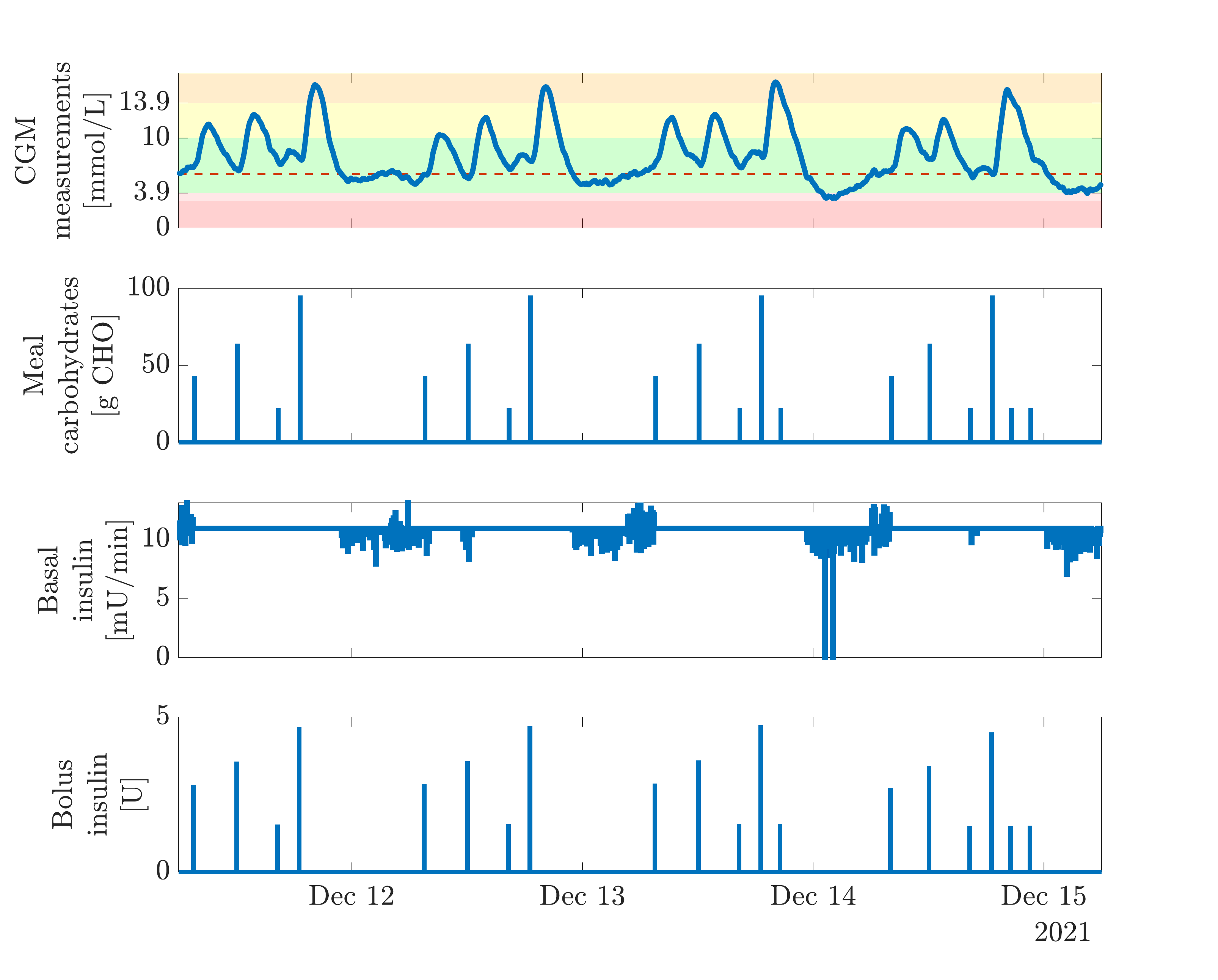}
	\caption{CGM measurements (top), meal carbohydrate contents (second from the top), basal insulin infusion rate (third from the top), and insulin boluses (bottom) for a single participant (represented with model~B) over a period of 4 days.}
	\label{fig:single:patient}
\end{figure}

We demonstrate the virtual clinical trial using a simple algorithm based on physiological insight and concepts from proportional-integral-derivative (PID) control. It uses two integrators (I-controllers) to estimate the basal and bolus insulin requirements, and a PD-controller to mitigate smaller variations in the blood glucose. It also uses deadbands, error truncation, a hypoglycemia amplification factor, switching logic, as well as a few heuristics. Fig.~\ref{fig:single:patient} shows how the controller doses insulin during 4~days of the virtual clinical trial for a single participant.

We evaluate the algorithm using the performance measures and targets described by~\citet{Holt:etal:2021b}. The measures are 1)~average glucose, 2)~glucose management index (GMI), 3)~glucose variation (GV) computed as the coefficient of variation, 3)~time above range (TAR), 4)~time in range (TIR), and 5)~time below range (TBR). All of these are based on CGM measurements. The 5~ranges used to compute TAR, TIR, and TBR are as follows (all values are in mmol/L, and the colors are used throughout this section). Level~2 hypoglycemia: $[0, 3[$ (red). Level~1 hypoglycemia: $[3, 3.9[$ (light red). Normoglycemia: $[3.9, 10]$ (green). Level 1 hyperglycemia, $]10, 13.9]$ (yellow). Level 2 hyperglycemia: $]13.9, \infty[$ (orange).

For model~A, $\tau_{IG}$ has the same value for all participants. All other parameters are sampled from normal distributions based on the sample means and variances of the parameter values presented by~\citet{Hovorka:etal:2002}. For model~B, we generate normal distributions based on the means and variances reported by~\citet{Colmegna:etal:2020}. For the remaining parameters, we use the means and variances of the values reported by~\citet{Kovatchev:etal:2010} and~\citet{DallaMan:etal:2007}. In both cases, we discard parameter sets that deviate more than one standard deviation from the mean or lead to a basal rate lower than 0.4~U/h (corresponding to a glucose concentration of 6~mmol/L). For model~A, we also require that all time constants (including the inverses of rate parameters) are within one order of magnitude of the mean.

\begin{figure*}[tb]
	\centering
	\includegraphics[trim=30 20 18 30, clip, width=0.1595\textwidth]{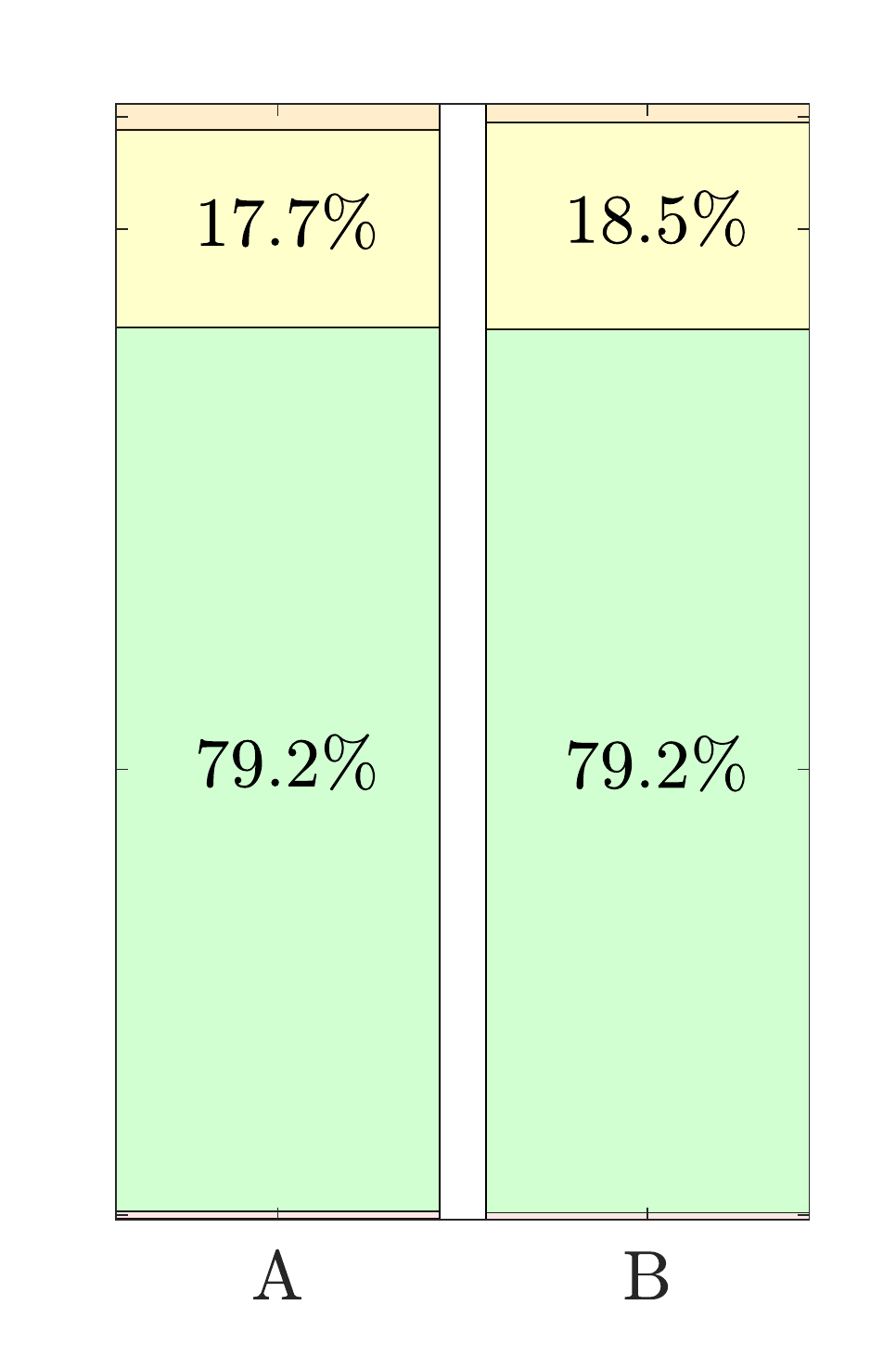}~
	\includegraphics[trim=30 20 16 30, clip, width=0.161\textwidth]{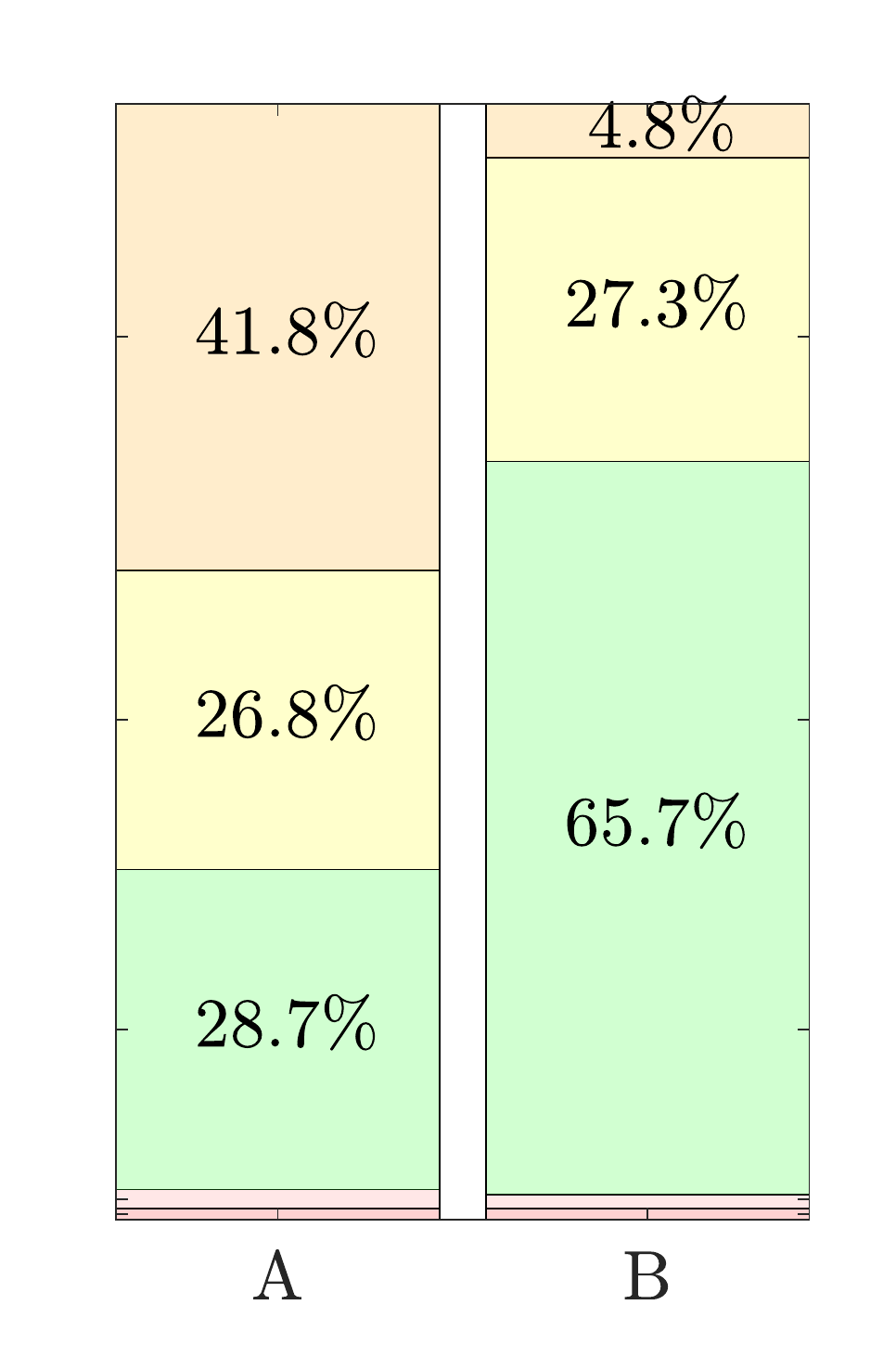}~
	\includegraphics[trim=70 20 80  5, clip, width=0.670\textwidth]{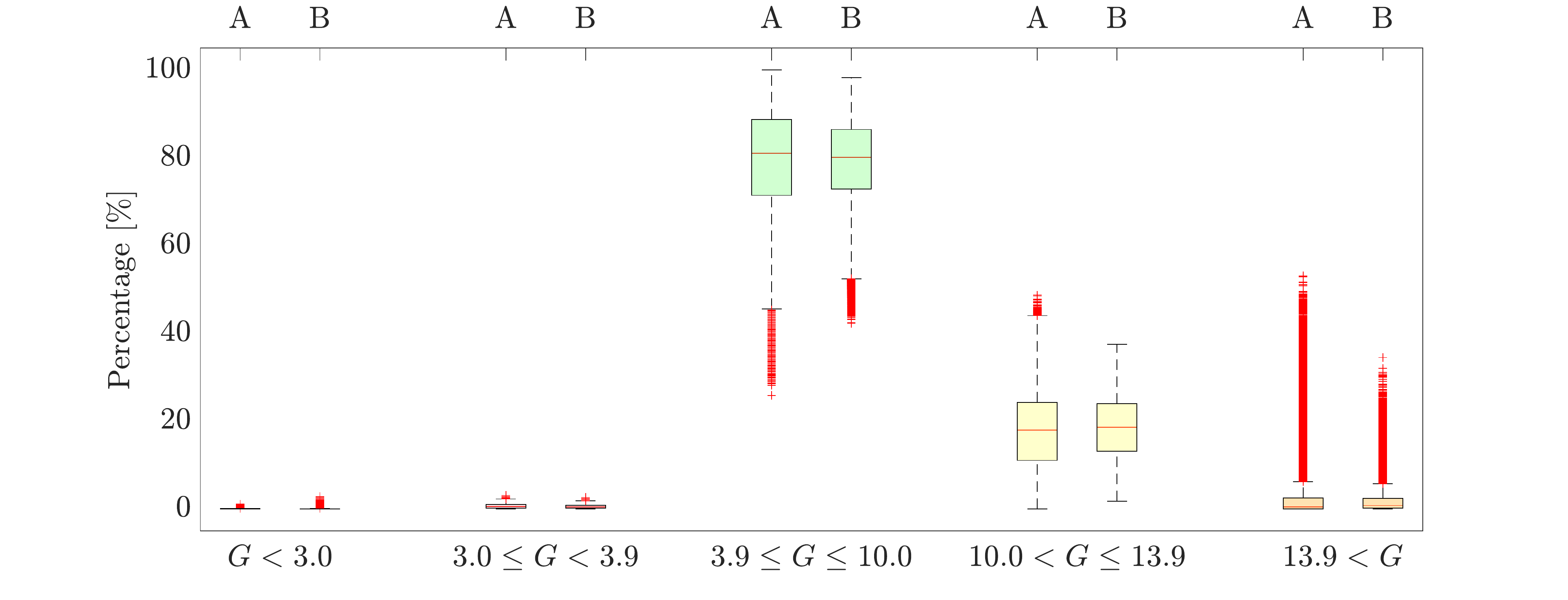}
	\caption{Left: Mean TIRs. Middle: TIRs for the worst-case participants. Right: Box plots of the TIRs for all participants.}
	\label{fig:tir}
\end{figure*}

\begin{figure}[tb]
	\centering
	\includegraphics[trim=60 35 80 45, clip, width=0.45\textwidth]{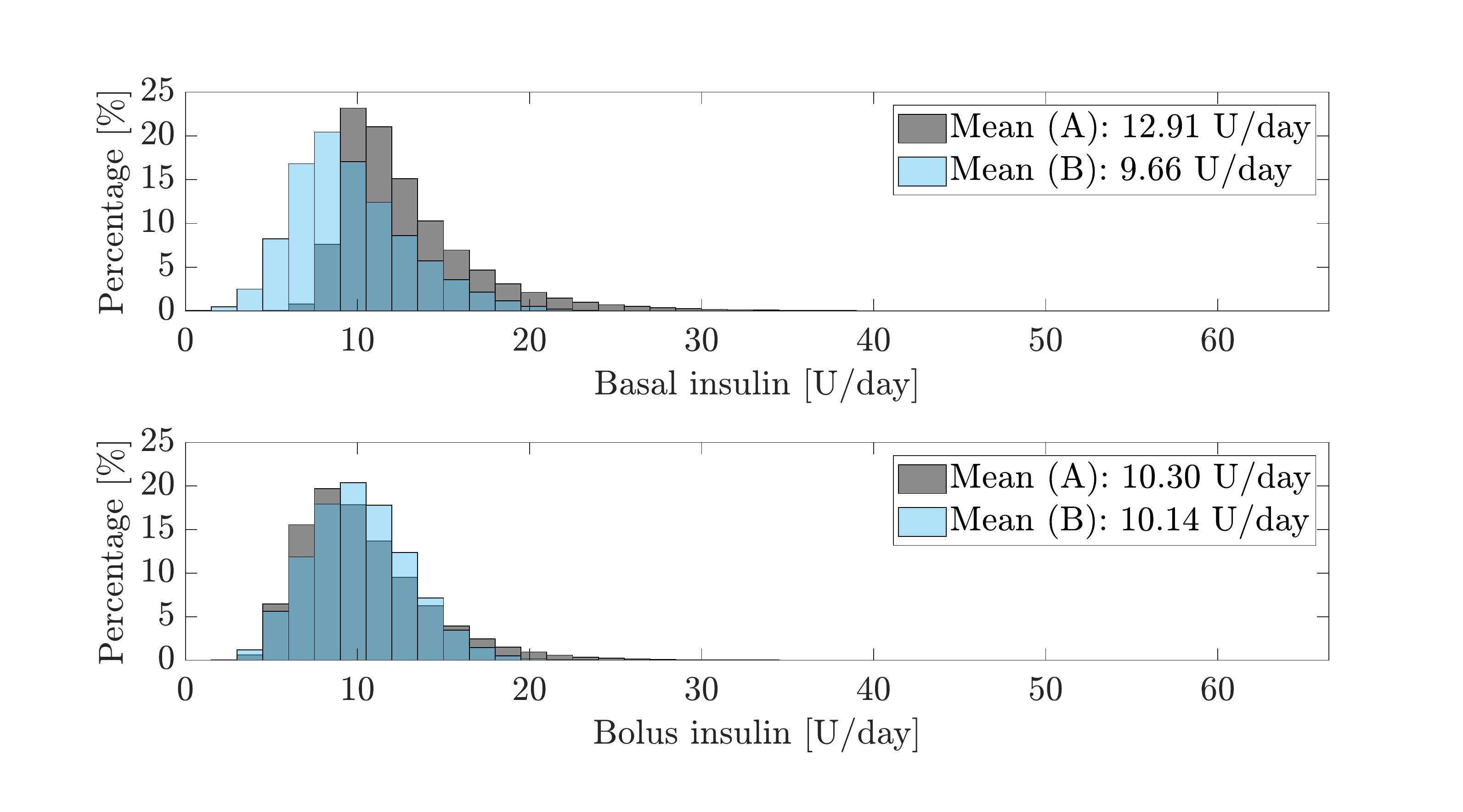}
	\caption{Histograms of the insulin basal and bolus TDDs.}
	\label{fig:tdd}
\end{figure}

\begin{figure}[tb]
	\centering
	\includegraphics[trim=50 5 75 20, clip, width=0.45\textwidth]{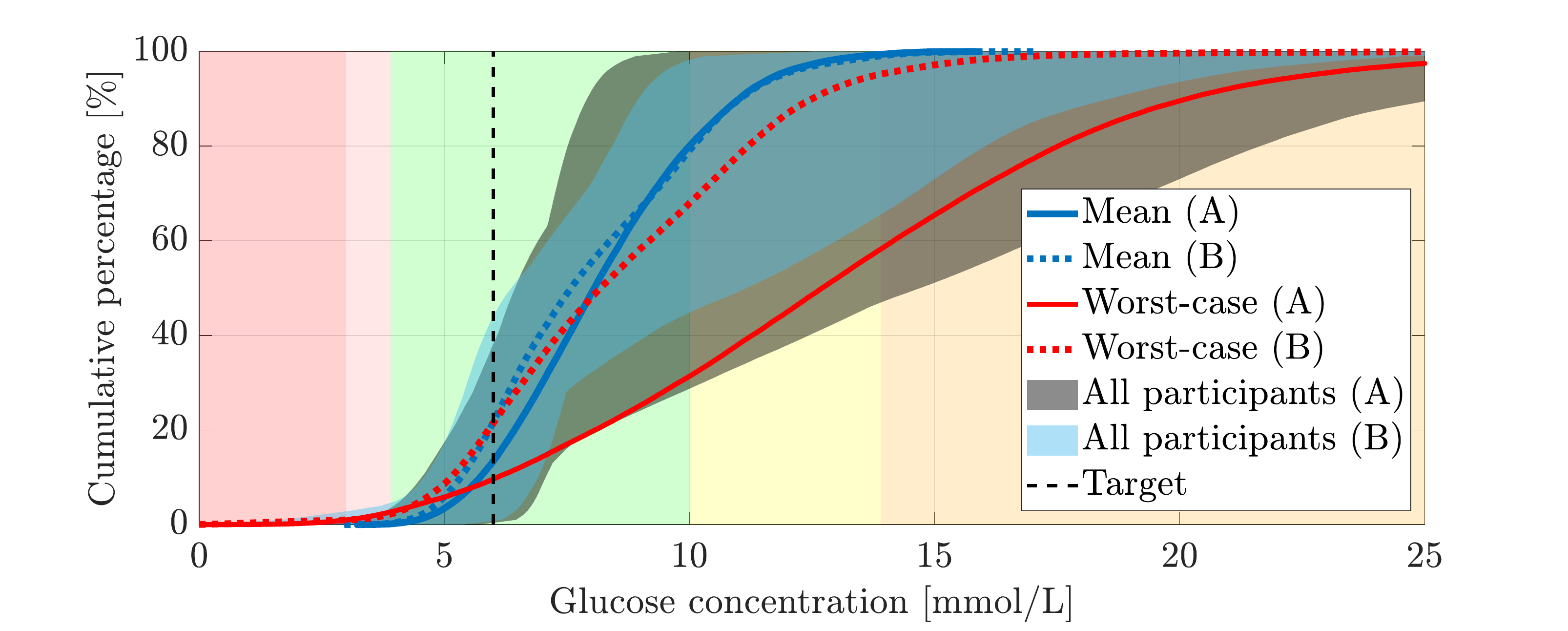}
	\caption{Cumulative distributions of the CGM measurements.}
	\label{fig:cumulative:distribution}
\end{figure}

Fig.~\ref{fig:tir} shows the mean and worst-case (i.e., the participant with the lowest CGM measurement) time in ranges (TIRs) as well as a box plot of the distributions of the TIRs. The mean TIRs are remarkably similar for the two models. The TIRs of the worst-case participants can be used to identify shortcomings. Here, the two worst-case participants appear to be problematic for different reasons. The box plot gives a more thorough overview and illustrates that, for model~B, the distributions of TIR and the times in level~1 and~2 hyperglycemia are more narrow. Fig.~\ref{fig:tdd} shows the total daily doses (TDDs) of basal and bolus insulin. For model~B, the participants get less basal insulin but almost the same amount of bolus insulin. Fig.~\ref{fig:cumulative:distribution} shows the cumulative distributions of the CGM measurements. It also indicates that the worst-case participants are qualitatively different. Finally, Table~\ref{tab:diabetes:treatment:targets} shows that the AP performs better for model~B. Presumably, the fact that the AP was tuned for model~A is outweighed by the (seemingly) smaller variations between the participants for model~B, which could be due to the parameter distributions.
\begin{table}[t]
	\centering
	\caption{Comparison based on the glycemic targets described by~\citet{Holt:etal:2021b}.}
	\label{tab:diabetes:treatment:targets}
	\begin{tabular}{llrr}
		Quantity 					& Target 		& \multicolumn{1}{c}{A} & \multicolumn{1}{c}{B} \\
		\hline \\[-9pt]
		Average glucose 			& $< 154$~mg/dL &  76.87\% 	&  91.09\%  \\
		GMI 						& $<     7\%$   &  77.48\% 	&  91.46\%  \\
		GV 							& $\leq 36\%$ 	&  91.37\% 	&  92.38\%  \\
		\hline \\[-9pt]
		TAR (level 2) 				& $<  5\%$ 		&  85.80\% 	&  92.55\%  \\
		TAR (level 1 and 2) 		& $< 25\%$ 		&  68.92\%  &  70.76\%  \\
		TIR\; (normoglycemia)		& $> 70\%$ 	    &  78.01\%  &  82.79\%  \\
		TBR (level 1 and 2) 		& $<  4\%$ 		& 100.00\%  & 100.00\%  \\
		TBR (level 2) 				& $<  1\%$ 		& 100.00\%  &  99.98\%  \\
		\hline \\[-9pt]
		\textbf{All targets} 		& 				&  68.10\% 	&  70.66\%  \\
		\hline
	\end{tabular}
\end{table}

    \section{Conclusions}\label{sec:conclusions}
We present an approach for conducting large-scale long-term virtual clinical trials for evaluating APs prior to the actual clinical trials. The participants are represented by multiple mathematical models consisting of stochastic differential equations. We use Monte Carlo closed-loop simulations, implemented with HPC software and hardware, to compute detailed treatment statistics. Finally, we present results from a virtual clinical trial with one~million participants (represented by two mathematical models) over 52 weeks, which is conducted in 2~h 9~min.

    % BIBLIOGRAPHY
    \bibliography{./ref/References}
\end{document}